\documentclass[11pt]{amsart}

\usepackage{amscd,amsmath,amssymb}
\usepackage[mathscr]{eucal}
\usepackage[all]{xy}
\usepackage[charter]{mathdesign}

\makeatletter \@addtoreset{equation}{section}\makeatother

\allowdisplaybreaks[1]

\newtheorem{theorem}{Theorem}[section]
\newtheorem{lemma}[theorem]{Lemma}
\newtheorem{proposition}[theorem]{Proposition}
\newtheorem{corollary}[theorem]{Corollary}

\newtheorem{remark}[theorem]{Remark}

\newcommand{\ZZ}{{\mathbb{Z}/2}}
\newcommand{\C}{{\mathbb{C}}}

\newcommand{\K}{{{k}}}
\newcommand{\End}{{\mathtt{End}}}
\newcommand{\lm}{{\mu}}
\newcommand{\ld}{{\delta}}
\newcommand{\bd}{{b(\delta)}}
\newcommand{\bm}{{b(\mu)}}
\newcommand{\un}{{\mathbf{1}\,}}

\newcommand{\uu}{{(\!(u)\!)}}
\newcommand{\uuu}{{[\![u]\!]}}
\newcommand{\rC}{{\mathsf{C}}}
\newcommand{\rCn}{{\overline{\mathsf{C}}}}
\newcommand{\rCp}{{\overline{\mathsf{C}}^{\Pi}}}

\newcommand{\str}{{\mathrm{str}}}
\newcommand{\A}{{\mathnormal{A}}}
\newcommand{\U}{{\mathnormal{e}(\delta)}}
\newcommand{\V}{{\mathnormal{E}(\delta)}}
\newcommand{\cI}{{\mathfrak{I}}}
\newcommand{\ck}{{\mathbf{sh}}}
\newcommand{\cK}{{\mathbf{Sh}}}
\newcommand{\cH}{{{H}}}
\newcommand{\sh}{{{sh}}}
\newcommand{\Sh}{{{Sh}}}
\newcommand{\cO}{{\mathscr{O}}}

\title[On a Hodge theoretic property of the K\"{u}nneth map]{On a Hodge theoretic property of the K\"{u}nneth map in periodic cyclic homology}

\author{Dmytro Shklyarov}
\address{Freiburg Institute for Advanced Studies (FRIAS) and Mathematics Institute, University of Freiburg,  Germany}
\email{dmytro.shklyarov@math.uni-freiburg.de}
\begin{document}
\begin{abstract}
We prove that the standard K\"{u}nneth map in periodic cyclic homology of differential $\ZZ$-graded algebras is compatible with a generalization of the Hodge filtration and explain how this result is related to various Thom-Sebastiani type theorems in singularity theory.
\end{abstract}

\maketitle

\baselineskip 1.5pc

\section{Introduction}
It is a classical fact in algebraic geometry that the K\"{u}nneth map in the cohomology of complex algebraic varieties respects Deligne's mixed Hodge structures. Our main result -- Theorem \ref{main} -- can be viewed as a first step towards generalizing this fact to the realm of non-commutative geometry. In our approach, ordinary spaces get replaced by differential $\ZZ$-graded (henceforth, `dg')  algebras over a fixed ground field $\K$ of characteristic 0 and the role of the classical cohomology is played by the periodic cyclic homology. The latter is well known to satisfy the `K\"{u}nneth property' \cite{Em,HJ,Kas,Lod,Pu} and is, in addition, anticipated to carry a Hodge-like structure \cite{KKP,Kon,KS}. What we show in the present paper is that the non-commutative K\"{u}nneth map respects one of the ingredients of the Hodge-like structure on the periodic cyclic homology. This ingredient, called the (formal) de Rham data \cite{KKP}, is an analog of the classical Hodge filtration.

To imagine what the de Rham data look like, recall that the periodic cyclic homology of a dg algebra is a super vector space over the field  $\K(\!(u)\!)$ of formal Laurent series in a variable $u$ which has the interpretation of a certain periodicity map \cite{GJ,HJ,Lod}. It will be convenient for us to think of such super vector spaces as super vector bundles over the punctured formal $u$-disk. (In order to stay in the familiar setting of finite rank bundles,  let us restrict ourselves to considering only those dg algebras whose periodic cyclic homology has finite dimension over $\K(\!(u)\!)$.)  For example, in these terms the aforementioned K\"{u}nneth property of the periodic cyclic homology is the claim that the tensor product of the bundles associated with two dg algebras is naturally isomorphic to the bundle associated with the tensor product of the dg algebras; an explicit canonical isomorphism, which we call {\it the} non-commutative K\"{u}nneth map, can be found in 
 \cite{Lod}.  The de Rham data  on the periodic cyclic homology of a dg algebra comprises two components: a canonical {\it connection} on the underlying bundle and a canonical {\it extension} of the bundle over the non-punctured formal $u$-disk (the extension is simply the image of the negative cyclic homology in the periodic one; we refer to \cite{KKP,Shk} for more details). That the non-commutative K\"{u}nneth map is compatible with the extensions is an obvious consequence of the explicit formula: roughly, the K\"{u}nneth map at the level of complexes is regular at $u=0$. The non-trivial part of our result is that the map is compatible with the connections.

While the compatibility of the classical K\"{u}nneth map with the mixed Hodge structures on the cohomology of varieties is a good piece of motivation for the subject of the present paper, there are  results in geometry of which our theorem is a {\it direct} generalization. Namely, our actual aim was to obtain abstract algebraic versions of various Thom-Sebastiani-type results in singularity theory such as, for instance, the Thom-Sebastiani formula for Steenbrink's Hodge filtration on the vanishing cohomology of isolated singularities \cite[Sect.8]{SS}. The latter calculates the Hodge filtration for the direct sum $f\oplus g$ of two singularities in terms of the same data for $f$ and $g$. This result, as well as the original Thom-Sebastiani theorem and some other facts of similar nature, can be deduced from a K\"unneth property for the Gauss-Manin systems and the Brieskorn lattices associated with the singularities (cf. Lemma 8.7 in {\it loc.cit.}) or, equivalently, for the Fourier-Laplace transforms thereof (cf. \cite[Sect. 3]{Sa}). This K\"unneth property is a special case of our result, as we will explain in Section \ref{mainthm}. 

{\bf Conventions.} In Section \ref{cycl}, we work over an arbitrary field $\K$ of characteristic 0. In Section \ref{mainthm} the ground field is $\C$. Given a $\ZZ$-graded space $V$, $\Pi V$ will stand for $V$ with the reversed $\ZZ$-grading, the parity of $v\in V$ will be denoted by $|v|$, and the corresponding element of $\Pi V$ will be denoted by $\Pi v$. We will follow all the standard conventions of super-linear algebra (such as the Koszul rule of signs, etc.). Finally, only unital dg algebras will be considered.

\medskip

{\bf Acknowledgement.} This research was supported by 
 the ERC Starting Independent Researcher Grant StG No. 204757-TQFT (K.~Wendland PI) and a Fellowship from the Freiburg Institute for Advanced Studies (FRIAS), Freiburg, Germany.

\section{Preliminaries and the main theorem}\label{cycl}
\subsection{Negative and periodic cyclic homology} Let us start by recalling the definition of the negative and periodic cyclic homology of a dg algebra (modulo minor details, our definitions agree with those in \cite{Get,GJ}).

Let $A=(A,d)$ stand for a dg algebra and set $\rC (A):=\bigoplus_{n\geq0} A\otimes (\Pi A)^{\otimes n}$. The latter is the underlying $\ZZ$-graded space of the Hochschild complex of $A$. The Hochschild differential, $b$, can be conveniently written in terms of some auxiliary operator on $\rC (A)$ which we will now introduce.

Writing the elements of $A\otimes (\Pi A)^{\otimes n}$ as $a_0[a_1|a_2|\ldots |a_n]$, we set
\begin{eqnarray*}
\tau(a_0[a_1|a_2|\ldots |a_n])=(-1)^{|\Pi a_0|\sum_{i=1}^n|\Pi a_i|}a_1[a_2|\ldots |a_n|a_0]
\end{eqnarray*}
\begin{eqnarray*}
\ld^{(0)}(a_0[a_1|a_2|\ldots |a_n])=da_0[a_1|a_2|\ldots |a_n]
\end{eqnarray*}
\begin{eqnarray*}\lm^{(0)}(a_0[a_1|a_2|\ldots |a_n])=
\begin{cases} 
0 & n=0\\
(-1)^{|a_0|}a_0a_1[a_2|\ldots |a_n] & n\geq1\\
\end{cases}
\end{eqnarray*}
\begin{eqnarray*}
\ld^{(i)}:=\tau^{-i}\,\ld^{(0)}\, \tau^{i}, \quad \lm^{(i)}:=\tau^{-i}\,\lm^{(0)}\, \tau^{i},\quad i=1,
\ldots,n
\end{eqnarray*}
Then the Hochschild differential $b$ is given by
$b=\bd+\bm$
where 
\[\bd=\sum_{i=0}^{n}\ld^{(i)},\quad\bm=\sum_{i=0}^{n}\lm^{(i)}\]

The {\it negative} cyclic homology of $A$ is the cohomology of the complex $(\rC (A)\uuu, b+uB)$ where $u$ is an even variable, $\rC (A)\uuu$ is the space of formal series with coefficients in $\rC (A)$, and $B=(1-\tau^{-1})hN$ 
where 
\[N=\sum_{i=0}^{n}\tau^{i},\quad h(a_0[a_1|a_2|\ldots |a_n])=1[a_0|a_1|a_2|\ldots |a_n]\]
The {\it periodic} cyclic homology is defined similarly, with the formal series in $u$ replaced by the formal Laurent series.

We will be working mostly with the {\it normalized} version of the complex. Namely, consider the subspace $\rC^{\prime}(A)\subset\rC(A)$ spanned by $a_0[a_1|\ldots |a_n]$ with at least one of $a_1,\ldots,a_n$ equal to 1. This subspace is preserved by both $b$ and $B$. Let 
\[
\rCn(A):=\rC(A)/\rC^{\prime}(A)=\bigoplus_{n\geq0} A\otimes (\Pi\, \overline{A})^{\otimes n}, \quad \overline{A}:=A/\K
\] 
We will use the same notation $b$, $B$ for the induced differentials on $\rCn(A)$ (the formula for $B$ on the quotient simplifies to
$B=hN$). The natural morphism of complexes
\[(\rC (A)\uuu, b+uB)\to(\rCn(A)\uuu, b+uB),\]
is known to be a quasi-isomorphism.

\subsection{The non-commutative K\"{u}nneth map}
Let us now recall the explicit formula for the non-commutative K\"{u}nneth map from\footnote{Formally speaking, the dg case is not discussed in \cite{Lod}. However, once the non-trivial grading is taken care of (cf. \cite{GJ,Ts}), the formula works in the dg setting as well.} \cite[Sect. 4.3]{Lod}. 

Consider two dg algebras, $A'$ and $A''$. Define 
\[\ck: \rCn(A')\otimes \rCn(A'')\to\rCn(A'\otimes A'')\]
as follows. For $a'_0[a'_1|\ldots|a'_n]\in \rCn(A')$ and $a''_0[a''_1|\ldots |a''_m]\in \rCn(A'')$ set
\begin{multline}\label{ku}
\ck(a'_0[a'_1|\ldots |a'_n]\otimes
a''_0[a''_1|\ldots|a''_m])\\=(-1)^{\ast}(a'_0\otimes a''_0)\sh[a'_1\otimes1|\ldots|a'_n\otimes1|1\otimes a''_1|\ldots|1\otimes a''_m]
\end{multline}
Here $\ast=|a''_0|(|\Pi a'_1|+\ldots+|\Pi a'_n|)$ and $\sh$ stands for the sum (with signs) over all the $(n,m)$-shuffles, i.e. the permutations that shuffle the $a'$-terms with the $a''$-terms while preserving the order of the former and the latter. The signs are computed by the rule that the transposition
$[\,\ldots|x|y|\ldots\,]\to[\,\ldots|y|x|\ldots\,]$ contributes $(-1)^{(|x|+1)(|y|+1)}$ (for $x$ and $y$ having definite parities).

The key property of $\ck$ is that it commutes with the Hochschild differentials 
\begin{equation*} b\,\ck=\ck\,(b\otimes 1+1\otimes b)\end{equation*}
and induces a quasi-isomorphism of the corresponding complexes (this is the `K\"{u}nneth theorem for Hochschild homology'). However, it does not induce a morphism of the {\it cyclic} complexes since it does not respect $B$. A first-order correction (in $u$) to $\ck$ turns out to solve this problem. Namely, let
\begin{multline*}
\cK(a'_0[a'_1|\ldots |a'_n]\otimes
a''_0[a''_1|\ldots|a''_m])\\=(-1)^{\ast\ast}(1\otimes 1)\Sh[a'_0\otimes1|\ldots|a'_n\otimes1|1\otimes a''_0|\ldots|1\otimes a''_m]
\end{multline*}
with $\ast\ast=|a'_0|+|\Pi a'_1|+\ldots+|\Pi a'_n|$ and the operator $\Sh$ defined by the same formula as $\sh$ but with the sum extended over all the {\it cyclic} shuffles. The latter cyclically permute $a'_0, \ldots,a'_n$ and $a''_0, \ldots,a''_m$, and then shuffle the $a'$-terms with the $a''$-terms so that $a'_0$ stays to the left of $a''_0$ (see \cite[Sect. 4.3.2]{Lod} and \cite[Sect. 4]{GJ}). The operator $\ck+u\cK$ commutes with the cyclic differentials
\begin{equation*} (b+uB)\,(\ck+u\cK)=(\ck+u\cK)\,((b+uB)\otimes 1+1\otimes (b+uB))\end{equation*}
and we obtain a morphism 
\[\left((\rCn(A')\otimes\rCn(A''))\uuu, (b+uB)\otimes 1+1\otimes (b+uB)\right)\to(\rCn(A'\otimes A'')\uuu, b+uB)\]
In general, it is only after passing to the formal Laurent series that it becomes a quasi-isomorphism. This is the `K\"{u}nneth theorem for periodic cyclic homology'.

\subsection{The canonical connection on the periodic cyclic homology}
The next definition we need to state our main result is that of the canonical connection on the periodic cyclic homology. 

We will need some terminology from \cite{Shk}. Namely, let $(\rC, b, B)$ be a mixed complex, i.e. a $\ZZ$-graded vector space endowed with an (ordered) pair of anti-commuting differentials. Then a {\it $u$-connection} on this mixed complex is an (even) differential operator of the form $\nabla=\frac{d}{du}+\A(u)$, $\A(u)\in\mathrm{End}(\rC)(\!(u)\!)$ satisfying
\begin{equation}\label{u-conn}
[\nabla\,,\, b+uB\,]=\frac1{2u}(b+uB)
\end{equation}
In particular, a $u$-connection induces a connection on the cohomology of $(\rC(\!(u)\!), b+uB)$ viewed as a bundle over the formal punctured disk.

Let $A$ be a dg algebra. Consider the operators on $\rCn(A)$ whose restrictions to $A\otimes (\Pi\, \overline{A})^{\otimes n}$ are given by 
\begin{eqnarray*} \gamma = n\cdot id,\quad
\U=-\lm^{(0)}\ld^{(1)}, \quad
\V=-\sum_{i=1}^{n}\sum_{j=0}^{n-i}h\tau^{-j}\ld^{(i)}
\end{eqnarray*}

\begin{proposition}\label{pro}
\begin{eqnarray}\label{1}
\nabla_A=\frac{d}{du}+\frac{\U}{2u^2}+\frac{\V-\gamma}{2u}
\end{eqnarray}
is a $u$-connection on $(\rCn(A), b, B)$. 
\end{proposition}
The induced connection on the periodic cyclic homology is the canonical connection.

Proposition \ref{pro} is a consequence of \cite[Cor. 3.7]{Shk} where a $u$-connection $\nabla^{un}$ on the unnormalized cyclic complex is presented. Inspecting the formula for $\nabla^{un}$ shows that $\nabla^{un}$ preserves the subspace $\rC'(A)\subset\rC(A)$ and the induced operator on the quotient $\rCn(A)(\!(u)\!)=\rC(A)/\rC'(A)(\!(u)\!)$ is precisely (\ref{1}). 

\begin{remark} {\rm The origin of the formula (\ref{1}), as well as other similar formulas in \cite{Shk}, is the non-commutative Cartan homotopy formula \cite{Get} (or, rather, a  special case of it known as Rinehart's formula). Namely, the relation (\ref{u-conn}) for the operator (\ref{1}) is equivalent to 
\begin{equation}\label{R}
[\U+u\V, b+uB]=u\bd
\end{equation}
which can be deduced from \cite[Eq.(2.1)]{Get} by substituting the differential $d$ of the dg algebra $A$ for the Hochschild cochain $D$ (in the notation of {\it loc.cit.}).}
\end{remark}

\subsection{The main result} As we explained in the Introduction, our main result says that the non-commutative K\"unneth map is compatible with the canonical connections. The actual statement is at the level of the cyclic complexes, and to formulate it we will need yet another (the last one) piece of terminology from \cite{Shk}. Namely, given two mixed complexes with $u$-connections, $(\rC', b', B', \frac{d}{du}+\A'(u))$ and $(\rC'', b'', B'', \frac{d}{du}+\A''(u))$, a {\it morphism} from the former to the latter is a $\K(\!(u)\!)$-linear morphism of complexes 
\[\psi(u): (\rC'(\!(u)\!), b'+uB')\to (\rC''(\!(u)\!), b''+uB'')\]
satisfying the following condition: the morphism
\[\frac{d\psi(u)}{du}+\A''(u)\psi(u)-\psi(u)\A'(u):(\rC'(\!(u)\!), b'+uB')\to (\rC''(\!(u)\!), b''+uB'')\]
is homotopic to 0.  This condition guarantees that the map, induced by $\psi(u)$ on the cohomology, respects the induced connections. 

Now we are able to state the main result:
\begin{theorem}\label{main}
The cyclic K\"unneth map $\ck+u\cK$ is a morphism from
\[\left(\rCn(A')\otimes\rCn(A''), b\otimes1+1\otimes b, B\otimes1+1\otimes B, \nabla_{A'}\otimes1+1\otimes \nabla_{A''}\right)\]
to
$(\rCn(A'\otimes A''), b,B,\nabla_{A'\otimes A''})$.
\end{theorem}
\begin{remark}\label{chc} {\rm Before proving the theorem, we would like to point out that all the constructions and conclusions in this section remain valid if we replace the ordinary cyclic mixed complexes $(\rCn(A), b, B)$ with their completed versions $(\rCp(A), b, B)$ where
\[\rCp(A):=\prod_{n\geq0} (A\otimes (\Pi\, \overline{A})^{\otimes n})^{\rm even}\bigoplus \prod_{n\geq0} (A\otimes (\Pi\, \overline{A})^{\otimes n})^{\rm odd}\]
Moreover, all the constructions are compatible with the canonical morphism 
\begin{equation}\label{compl}(\rCn(A)\uuu, b+uB)\to(\rCp(A)\uuu, b+u B)\end{equation}
Note that, in general, (\ref{compl}) is not a quasi-isomorphism (although it is a quasi-isomorphism in a number of cases of interest; cf. \cite{PP}).}
\end{remark}
\noindent{\bf Proof of Theorem \ref{main}.}
Let us write expressions like 
\[b\,\ck-\ck\,(b\otimes 1+1\otimes b),\quad \U\,\cK-\cK\,(\U\otimes 1+1\otimes\U),\quad \text{etc.}\]
using the commutator notation:
$[b,\ck]$, $[\U,\cK]$ etc. (the `commutator' here stands for the {\it super}-commutator when both operators involved are odd).

We have to prove that the morphism
\[2u^2\frac{d(\ck+u\cK)}{du}+[\U+u(\V-\gamma),\ck+u\cK]\]
is homotopic to 0. Consider its $u$-expansion:
\begin{eqnarray*}
[\U,\ck]+u\left([\U,\cK]+[\V-\gamma,\ck])\right)+u^2\left([\V-\gamma,\cK]+2\cK\right)
\end{eqnarray*}
Obviously, the coefficient at $u^2$ equals 0. Indeed, $[\gamma,\cK]=2\cK$ and $[\V,\cK]=0$
simply by the definition of the operators involved and the fact that we are working with the normalized complexes. 
The coefficient at $u^0$ is also easily seen to vanish. The proof of this is based on the explicit formula for $\U$ 
\[\U(a_0[a_1|\ldots|a_n])=a_0da_1[a_2|\ldots|a_n]\]
and the obvious fact that the sum in the right-hand side of (\ref{ku}) involves elements of the following two types only:
\[\pm (a'_0\otimes a''_0)[a'_1\otimes 1|\ldots]\quad \text{or}\quad \pm (a'_0\otimes a''_0)[1\otimes a''_1|\ldots]
\]
Finally, note that the coefficient at $u$ equals
\begin{equation}\label{d}
[\U,\cK]+[\V,\ck]
\end{equation}
since, clearly, $[\gamma, \ck]=0$.
Thus, we need to show that the operator (\ref{d}), viewed as a morphism 
\[\left((\rCn(A')\otimes\rCn(A''))(\!(u)\!), (b+uB)\otimes 1+1\otimes (b+uB)\right)\to(\rCn(A'\otimes A'')(\!(u)\!), b+uB),\]
is homotopic to 0.

Let  
$
\cH:\rCn(A')\otimes\rCn(A'')\to\rCn(A'\otimes A'')
$
denote the odd operator whose restriction to $\left(A'\otimes (\Pi\, \overline{A'})^{\otimes n}\right)\otimes \left(A''\otimes (\Pi\, \overline{A''})^{\otimes m}\right)$
is given by 
\[
\cH=\sum_{i=1}^{n}\sum_{r=0}^{n-i}\sum_{s=0}^{m}\ck^{(r,s)}(\tau^{-r}\ld^{(i)}\otimes \tau^{-s})+\sum_{i=1}^{m}\sum_{r=0}^{n}\sum_{s=0}^{m-i}\ck^{(r,s)}(\tau^{-r}\otimes \tau^{-s}\ld^{(i)})
\]
where
\begin{multline*}
\ck^{(r,s)}(a'_0[a'_1|\ldots |a'_n]\otimes
a''_0[a''_1|\ldots|a''_m])=\\(-1)^{\ast\ast}(1\otimes 1)\sh^{(r,s)}[a'_0\otimes 1|\ldots|a'_n\otimes 1|1\otimes a''_0|\ldots|1\otimes a''_m]
\end{multline*}
with $\ast\ast$ being the same as in the definition of $\cK$ and $\sh^{(r,s)}$ denoting the sum over those (non-cyclic) $(n+1,m+1)$-shuffles that do not switch $a'_r\otimes1$ and $1\otimes a''_s$. 

It is not hard to prove that $[b+uB,\cH]=[\bm,\cH]$. Indeed, $[\bd,\cH]=0$ since $\bd$ commutes with the shuffles and the cyclic permutations, and anti-commutes with all the $\ld^{(i)}$. Also, $[B,\cH]=0$ since we are working with the normalized complexes. On the contrary, the proof of the following lemma is a rather long and tedious calculation; an interested reader may find it in Appendix \ref{mainap}. 

\begin{lemma}\label{l1}
$[\U,\cK]+[\V,\ck]=
[\bm,\cH]+\ck(\bd\otimes B) 
$
\end{lemma}

By combining the above claims, we obtain
\[
[\U,\cK]+[\V,\ck]=
(b+uB)\cH+\cH((b+uB)\otimes1+1\otimes(b+uB))+\ck(\bd\otimes B) 
\] 
Thus, it remains to show that the operator $\ck(\bd\otimes B)$ is homotopic to 0. Since we are working with the normalized complexes,
\[
\ck(\bd\otimes B)=(\ck+u\cK)(\bd\otimes B)=(\ck+u\cK)(\bd\otimes1)(1\otimes B)
\]
Observe that the operators $\ck+u\cK$ and $1\otimes B$ (anti-)commute with the cyclic differential while $\bd\otimes1$ can be written as the commutator $[b+uB, H']$ (see (\ref{R})).

\section{Relation to Thom-Sebastiani type theorems in singularity theory}\label{mainthm}
From now on, the ground field is $\C$.
\subsection{The Gauss-Manin connections and Thom-Sebastiani type results}
Let $f=f(x_1,\ldots,x_k)$ be a polynomial having an isolated critical point at the origin, with the corresponding critical value equal to $0$. Associated with the germ of $f$ at the origin is the  Gauss-Manin system -- a holonomic $D$-module in one variable -- which, together with the so-called Brieskorn lattice, encodes most of the Hodge theoretic invariants of the critical point. 
We refer the reader to \cite{SS} for a detailed exposition; a streamlined reminder can be found in \cite[Sect. 5]{Shk}. Essentially the same amount of information is contained in the (formal) Fourier-Laplace transforms of the Gauss-Manin system and the Brieskorn lattice, and it is these invariants that we will be interested in in this section. Let us recall their explicit description.

The starting point is the (formal) {\it twisted de Rham cohomology} of the critical point, i.e. the cohomology of the complex $({\Omega}^{an}_{\C^k, 0}\uu, -df+ud)$. Here ${\Omega}^{an}_{\C^k, 0}$ denotes the $\ZZ$-graded space of germs at the origin $0\in\C^k$ of analytic differential forms. The cohomology $H^{\ast}({\Omega}^{an}_{\C^k, 0}\uu, -df+ud)$ is known to be non-trivial in degree $k$ only. The Fourier-Laplace transformed Gauss-Manin system is the $\C\uu$-linear space $H^{k}({\Omega}^{an}_{\C^k, 0}\uu, -df+ud)$ equipped with a connection $\nabla_f^{GM}$ induced by the differential operator $\frac{d}{du}+\frac{f}{u^2}$ on ${\Omega}^{an}_{\C^k, 0}\uu$. The Brieskorn lattice transforms into the $\C\uuu$-lattice $H^{k}({\Omega}^{an}_{\C^k, 0}\uuu, -df+ud)\subset H^{k}({\Omega}^{an}_{\C^k, 0}\uu, -df+ud)$.

Let $g=g(y_1,\ldots,y_l)$ be another polynomial with the same property, i.e. having an isolated critical point at the origin, with the corresponding critical value equal to $0$. Recall that $f\oplus g$ stands for the polynomial $f(x_1,\ldots,x_k)+g(y_1,\ldots,y_l)$ on $\C^{k}\times\C^{l}$ (clearly, it also has an isolated singularity at the origin). It is not hard to see that the map
\[
\bigwedge:{\Omega}^{an}_{\C^k, 0}\otimes{\Omega}^{an}_{\C^l, 0}
\to{\Omega}^{an}_{\C^k\times\C^l, 0},\quad \omega'(x_1,\ldots,x_k)\otimes \omega''(y_1,\ldots,y_l)\mapsto \omega'(x_1,\ldots,x_k)\wedge\omega''(y_1,\ldots,y_l)
\]
induces an isomorphism 
\[
H^{k}({\Omega}^{an}_{\C^k, 0}\uuu, -df+ud)\otimes_{\C\uuu} H^{l}({\Omega}^{an}_{\C^l, 0}\uuu, -dg+ud)\to H^{k+l}({\Omega}^{an}_{\C^k\times\C^l,0}\uuu, -d(f\oplus g)+ud)
\]
which transforms $\nabla_f^{GM}\otimes 1+1\otimes \nabla_g^{GM}$ into $\nabla_{f\oplus g}^{GM}$. 

As explained in \cite[Sect. 8]{SS}, the above `Thom-Sebastiani theorem' for the Brieskorn lattices and the Gauss-Manin connections implies a number of other Thom-Sebastiani type formulas. The aim of this section is to explain that Theorem \ref{main} contains the Thom-Sebastiani theorem for the Gauss-Manin connections as a special case.

\subsection{Reminder on dg algebras associated with isolated critical points} In this section, we will switch to the algebraic setting and, in particular, will view polynomials as germs of algebraic functions: $f\in \cO_{\C^k,0}$ where $\cO_{\C^k,0}$ stands for the local algebra of the affine variety $\C^k$ at $0$. Accordingly, we will work with germs of algebraic differential forms which will be denoted simply by ${\Omega}_{\C^k, 0}$. Note that the formal twisted de Rham cohomology does not see the difference between the analytic and the algebraic forms: as explained in \cite{Sch}, the natural maps
\[
H^{k}({\Omega}_{\C^k, 0}\uuu, -df+ud)\to H^{k}({\Omega}^{an}_{\C^k, 0}\uuu, -df+ud)\to H^{k}({\Omega}^{form}_{\C^k, 0}\uuu, -df+ud)
\]
are quasi-isomorphisms (here ${\Omega}^{form}_{\C^k, 0}$ stands for the space of formal germs of differential forms\footnote{In fact, we could have chosen to work in the purely formal setting, i.e. to view $f$ as an element of the algebra of formal series. The results to be discussed  hold true in this setting.}).

According to \cite{Dyck}, associated\footnote{The dg algebra is non-canonical but its dg Morita equivalence class coincides with that of the dg category of matrix factorizations of $f$ and, thus, is functorial in $f$. } with $f$ is a dg algebra $A_f$.  The underlying graded algebra is 
$
\cO_{\C^k,0}\otimes{\End} P_k
$ 
where $P_k$ is the space of polynomials in $k$ odd variables $\theta_1,\ldots,\theta_k$. To equip it with a differential, we need to pick a decomposition 
\begin{equation}\label{dec}
f=x_1f_1+\ldots+x_kf_k,  \quad f_i\in \cO_{\C^k,0}.
\end{equation} 
Then the differential is the super-commutator with $D_f:=\sum_i\left(f_i\otimes\theta_i+x_i\otimes\partial_{\theta_i}\right)\in \cO_{\C^k,0}\otimes{\End} P_k$. 

The importance of this dg algebra for our purposes stems from two results obtained in \cite{Seg} and \cite{Shk}, respectively. Namely, in \cite{Seg} an explicit quasi-isomorphism 
\begin{equation}\label{qis1}
 \cI_f: (\rCn(A_f)\uuu, b+uB)\to({\Omega}_{\C^k,0}\uuu, -df+ud)
\end{equation}
was constructed. As a map of graded spaces, $\cI_f$ is the $\C\uuu$-linear extension of the composition 
\begin{equation*}\label{qis}
\rCn(\cO_{\C^k,0}\otimes{\End} P_k)\xrightarrow{\mathrm{exp}(-b(D_f))}\rCp(\cO_{\C^k,0}\otimes{\End} P_k)\stackrel{\str}\longrightarrow\rCp(\cO_{\C^k,0})\stackrel{\epsilon}\longrightarrow{\Omega}_{\C^k,0}
\end{equation*}
where 
\begin{itemize}
\item $\rCp$ is the completed version of $\rCn$ that we mentioned in Remark \ref{chc};\item $
b(D_f)(a_0[a_1|\ldots |a_n]):=\sum_{i=1}^{n+1}a_0[a_1|\ldots|a_{i-1}|D_f|a_{i}| \ldots|a_n]; 
$
\item 
$
\str (\phi_0\otimes T_0[\phi_1\otimes T_1|\ldots |\phi_n\otimes T_n]):= (-1)^{\sum_{i\,{\rm odd}}|T_i|}\str(T_0T_1\ldots T_n)\phi_0[\phi_1|\ldots |\phi_n] 
$\\ 
where $\phi_i\in \cO_{\C^k,0}$, $T_i\in{\End} P_k$;
\item
 $\epsilon$ is the Hochschild-Kostant-Rosenberg map: 
\[
\epsilon(\phi_0[\phi_1|\ldots |\phi_n])=\frac{1}{n!}\phi_0d\phi_1\wedge\ldots \wedge d\phi_n.
\]
\end{itemize}

In \cite{Shk}  $\cI_f$ was shown to define a morphism of mixed complexes with $u$-connections 
\[
\left(\rCn(A_f), b, B, \nabla_{A_f}\right)\to ({\Omega}_{\C^k,0}, -df, d, \frac{d}{du}+\frac{f}{u^2}-\frac{\gamma}{2u}),\quad \gamma|_{{\Omega}_{\C^k,0}^d}=d\cdot id
\] 
Let us denote by  $\nabla_f$ the connection on $H^{k}({\Omega}_{\C^k,0}\uu, -df+ud)$ induced by the operator
$
\frac{d}{du}+\frac{f}{u^2}-\frac{\gamma}{2u}.
$
Observe that $\nabla_f$ and $\nabla_f^{GM}$ differ only by a `Tate twist' $\frac{k}{2u}$. Thus, the above result means that, up to the Tate twist, $\nabla_f^{GM}$ can be recovered from the dg algebra associated with $f$.

Note  that the Tate twist is additive with respect to the direct sum of polynomials. Hence the Thom-Sebastiani theorem for $\nabla^{GM}_f$ and $\nabla^{GM}_g$ implies that for $\nabla_f$ and $\nabla_g$, and vice versa.

\subsection{Comparing the commutative and the non-commutative Thom-Sebastiani}
Let, as before, $f$ and $g$ be two polynomials with isolated critical points. 

Let us fix decompositions of the form (\ref{dec}) for $f$ and $g$. Then we automatically get a decomposition for ${f\oplus g}$, together with an obvious embedding of dg algebras $\iota:A_{f}\otimes A_{g}\to A_{f\oplus g}$. 

\begin{theorem}\label{TS}
The diagram of morphisms of complexes
\[
\xymatrix{ (\rCn(A_f)\otimes \rCn(A_g))\uuu, (b+uB)\otimes 1+1\otimes (b+uB)
\ar[r]^<<<<<<<<{\ck+u\cK} 
\ar[dd]_{\cI_f\otimes\cI_g}
&
\rCn(A_{f}\otimes A_{g})\uuu, b+uB \ar[d]_{\rCn(\iota)}\\
&
\rCn(A_{f\oplus g})\uuu, b+uB \ar[d]_{\cI_{f\oplus g}}\\
({\Omega}_{\C^k, 0}\otimes{\Omega}_{\C^l, 0})\uuu, (-df+ud)\otimes 1+1\otimes (-dg+ud)\ar[r]^<<<<<{\bigwedge}
&
{{\Omega}_{\C^k\times\C^l, 0}\uuu}, -d(f\oplus g)+ud}
\]
becomes commutative after passing to cohomology.
\end{theorem}
Since all the morphisms in the diagram preserve the $u$-connections, by localizing at $u$ we obtain
\begin{corollary}
 The Thom-Sebastiani theorem for the connections $\nabla_f$ and $\nabla_g$ -- and, consequently, for $\nabla^{GM}_f$ and $\nabla^{GM}_g$ -- is a special case of Theorem \ref{main}.
\end{corollary}
 
\noindent{\bf Proof of Theorem \ref{TS}.} As is shown in \cite{PP}, the canonical morphism (\ref{compl}) is a quasi-isomorphism for $A:=A_f$. Thus, in the above diagram, we may replace the ordinary cyclic complexes with their completed versions.
 
An easy calculation shows that 
\[b(D_f\otimes1+1\otimes D_g)\cdot\ck=\ck\cdot(b(D_f)\otimes1+1\otimes b(D_g)),\]
\[\str\cdot\ck=\ck\cdot(\str\otimes\str),\]
\[\epsilon\cdot\ck=\bigwedge\cdot\,(\epsilon\otimes\epsilon).\]
As a consequence, the following diagram (of maps of graded spaces) is commutative:
\[
\xymatrix{ \rCp(\cO_{\C^k,0}\otimes{\End} P_k)\otimes \rCp(\cO_{\C^l,0}\otimes{\End} P_l)
\ar[r]^<<<<{\ck}\ar[dd]_{\cI_f\otimes\cI_g}
&
\rCp(\cO_{\C^k,0}\otimes{\End} P_k\otimes \cO_{\C^l,0}\otimes{\End} P_l) \ar[d]_{\rCp(\iota)}\\
&
\rCp(\cO_{\C^k\times \C^l,0}\otimes{\End} (P_k\otimes P_l)) \ar[d]_{\cI_{f\oplus g}}\\
{\Omega}_{\C^k, 0}\otimes{\Omega}_{\C^l, 0}\ar[r]^<<<<<<<<<<<<<<<<<<<<{\bigwedge}
&
{{\Omega}_{\C^k\times\C^l, 0}}}
\]
It follows, in particular, that the composition $\cI_{f\oplus g}\cdot\rCp(\iota)\cdot\ck$ induces a morphism of complexes:
\[
\cI_{f\oplus g}\cdot\rCp(\iota)\cdot\ck\,((b+uB)\otimes 1+1\otimes (b+uB))=(-d(f\oplus g)+ud)\,\cI_{f\oplus g}\cdot\rCp(\iota)\cdot\ck
\]
Therefore, since $\cI_{f\oplus g}\cdot\rCp(\iota)\cdot(\ck+u\cK)$ commutes with the differentials as well, so does the composition $\cI_{f\oplus g}\cdot\rCp(\iota)\cdot\cK$. To prove the theorem, we need to show that the latter map is trivial at the level of cohomology.

Let us introduce the following decreasing filtration:
\[
F^p\rCp(A):= \prod_{n\geq p} (A\otimes (\Pi\, \overline{A})^{\otimes n})^{\rm even}\bigoplus \prod_{n\geq p} (A\otimes (\Pi\, \overline{A})^{\otimes n})^{\rm odd}
\]
Observe that, by the definitions of the operators involved,
\[
\cI_{f\oplus g}\cdot\rCp(\iota)\cdot\cK(F^p\rCp(A_f)\otimes F^q\rCp(A_g))\subset \cI_{f\oplus g}(F^{p+q+1}\rCp(A_{f\oplus g}))\subset\epsilon(F^{p+q+1}\rCp(\cO_{\C^k\times \C^l,0}))
\] 
which implies
\[
\cI_{f\oplus g}\cdot\rCp(\iota)\cdot\cK(F^p\rCp(A_f)\otimes F^q\rCp(A_g))=0,\quad p+q\geq k+l.
\]
Thus, to complete the proof, it would suffice to show that 
any class in $H^\ast(\rCp(A_f)\uuu, b+uB)$ (resp. $H^\ast(\rCp(A_g)\uuu, b+uB)$) can be represented by an element of $F^k\rCp(A_f)\uuu$ (resp. $F^l\rCp(A_g)\uuu$). 
To prove this, we need a bit more information about $\cI_f$. 

In fact, every map in the sequence
\[
\rCp(\cO_{\C^k,0}\otimes{\End} P_k)\uuu\xrightarrow{\mathrm{exp}(-b(D_f))}\rCp(\cO_{\C^k,0}\otimes{\End} P_k)\uuu\stackrel{\str}\longrightarrow\rCp(\cO_{\C^k,0})\uuu\stackrel{\epsilon}
\longrightarrow{\Omega}_{\C^k,0}\uuu
\]
is a morphism of mixed complexes (cf. \cite[Sect. 4.1]{Shk}), namely,
\[
\mathrm{exp}(-b(D_f)): (\rCp(A_f),\, b,\, B)\to (\rCp(\cO_{\C^k,0}\otimes{\End} P_k),\, \bm+b(f),\, B),
\]
\[
\str: (\rCp(\cO_{\C^k,0}\otimes{\End} P_k),\, \bm+b(f), \,B)\to (\rCp(\cO_{\C^k,0}),\, \bm+b(f), \,B),
\]
\[
\epsilon: (\rCp(\cO_{\C^k,0}), \,\bm+b(f), \,B)\to ({\Omega}_{\C^k,0},\, -df, \,d),
\]
where
\[
b(f)(a_0[a_1|\ldots |a_n]):=\sum_{i=1}^{n+1}(-1)^{\sum_{j=0}^{i-1}|\Pi a_j|}a_0[a_1|\ldots|a_{i-1}|f|a_{i}| \ldots|a_n].
\]
Moreover, each map induces a quasi-isomorphism at the level of the associated $\C\uuu$-linear complexes.

Now, let us pick a collection $f^{(1)},\ldots, f^{(\mu_f)}\in \cO_{\C^k,0}$ so that the corresponding classes in the Milnor algebra $\cO_{\C^k,0}/(\partial_{x_1}f,\ldots, \partial_{x_k}f)$ form a basis. Then the classes of the elements $\{f^{(j)}dx_1\wedge\ldots\wedge dx_k\}_{j=1,\ldots, \mu_f}$ form a basis of the free $\C\uuu$-module $H^{\ast}({\Omega}_{\C^k, 0}\uuu, -df+ud)$. Let
\[
{\bf f\,}^{(j)}=\sum_{\alpha=1}^{N_j} \phi^{(j)}_{0,\alpha}[\phi^{(j)}_{1,\alpha}|\ldots |\phi^{(j)}_{n,\alpha}],\quad  j=1,\ldots, \mu_f
\]
be any $\bm$-closed elements in $\rCp(\cO_{\C^k,0})$ such that $\epsilon({\bf f\,}^{(j)})=f^{(j)}dx_1\wedge\ldots\wedge dx_k$. Then, for any even idempotent $e\in {\End} P_k$ such that $\str(e)=1$, we obtain a collection of $\bm$-closed elements in $\rCp(\cO_{\C^k,0}\otimes {\End} P_k)$, namely,
\[
\dot{{\bf f}}\,^{(j)}=\sum_{\alpha=1}^{N_j} \phi^{(j)}_{0,\alpha}\otimes e[\phi^{(j)}_{1,\alpha}\otimes e|\ldots |\phi^{(j)}_{n,\alpha}\otimes e],\quad  j=1,\ldots, \mu_f.
\]
By (an easy extension of) the classical HKR theorem,  the map $\epsilon\cdot\str: (\rCp(\cO_{\C^k,0}\otimes{\End} P_k),\, \bm)\to ({\Omega}_{\C^k,0}, 0)$ is a quasi-isomorphism. In particular, any  $\bm$-closed element in $F^{p}\rCp(\cO_{\C^k,0}\otimes{\End} P_k)$, $p>k$, is $\bm$-exact. Using this observation, one can show that the elements $\dot{{\bf f}}\,^{(j)}$ extend to $(\bm+b(f)+uB)$-closed elements $\ddot{{\bf f}}\,^{(j)}$ in $F^{k}\rCp(\cO_{\C^k,0}\otimes{\End} P_k)\uuu$ such that $\ddot{{\bf f}}\,^{(j)}-\dot{{\bf f}}\,^{(j)}\in F^{k+1}\rCp(\cO_{\C^k,0}\otimes{\End} P_k)\uuu$. Since $\epsilon\cdot\str(\ddot{{\bf f}}\,^{(j)})=f^{(j)}dx_1\wedge\ldots\wedge dx_k$, the collection $\{\ddot{{\bf f}}\,^{(j)}\}_{j=1,\ldots, \mu_f}$ is a basis of the free $\C\uuu$-module $H^{\ast}(\rCp(\cO_{\C^k,0}\otimes{\End} P_k)\uuu,\bm+b(f)+uB)$. What remains is to observe that $\mathrm{exp}(b(D_f))$ preserves the filtration.

\appendix
\section{Proof of Lemma \ref{l1}}\label{mainap}
In this proof, all the operators will be viewed as operators on the {\it unnormalized} cyclic complexes. We will use the following shorthand notation: 
\begin{itemize}
\item Given an operator on $\rC (A)$ we equip the notation for the operator with the subscript ${n+1}$ to denote its restriction onto $A\otimes (\Pi A)^{\otimes n}$. 
\item We will omit the symbol $\otimes$ in the notation for the elements of $A'\otimes A''$. Namely, $a'$ (resp. $a''$) will stand  either for an element of $A'$ (resp. $A''$) or for the corresponding element $a'\otimes 1$ (resp. $1\otimes a''$) of $A'\otimes A''$, and $a'\otimes a''$ will be shortened to $a'a''$. 
\item The unit $1\otimes1$ of $A'\otimes A''$ will be denoted by $\un$. 
\end{itemize}

\begin{lemma}\label{l1ap}
\[
\V\ck-\cK(\U\otimes1+1\otimes\U)=
\bm^\dagger\cH+\cH(\bm\otimes1+1\otimes\bm)
\] 
where $\bm^\dagger$ stands for the following truncation of $\bm$: $\bm^\dagger_{n+1}=\sum_{i=1}^{n-1}\lm^{(i)}_{n+1}$.
\end{lemma}
\noindent{\bf Proof.} Let $\ck^{(r,s)}_0$ be the linear operator defined similarly to $\ck^{(r,s)}$ but with the sum taken only over those shuffles that put $a'_r$ and $a''_s$ next to each other: 
\[
\ck^{(r,s)}_0(a'_0[a'_1|\ldots |a'_n]\otimes a''_0[a''_1|\ldots|a''_m])=(-1)^{\ast\ast}(-1)^{\ast\ast\ast}\un[\ldots]
\]
where $\ast\ast$ is as in the definitions of $\cK$ and $\ck^{(r,s)}$, $[\ldots]$ stands for the following sum
\begin{equation}\label{k0} 
[\underbrace{a'_0|\ldots|a'_{r-1}|a''_{0}|\ldots|a''_{s-1}}_{\text{sum}\,
\text{over}\,(r,s)-\text{shuffles}}|a'_r|a''_s|\underbrace{a'_{r+1}|\ldots|a'_n|a''_{s+1}|\ldots|a''_m}_{\text{sum}\,
\text{over}\,(n-r,m-s)-\text{shuffles}}]
\end{equation}
and $(-1)^{\ast\ast\ast}$ is the sign coming from the shuffles that transform  $[a'_0|\ldots|a''_m]$ into (\ref{k0}). 
It is easy to see that
\begin{eqnarray*}
\bm^\dagger\ck^{(r,s)}=&&\lm^{(r+s+1)}\ck^{(r,s)}_0\\
&+&\ck^{(r-1,s)}\sum_{j=0}^{r-1}\lm_{n+1}^{(j)}\otimes1+\ck^{(r,s)}\sum_{j=r}^{n-1}\lm_{n+1}^{(j)}
\otimes1\\
&+&\ck^{(r,s-1)}\sum_{j=0}^{s-1}1\otimes \lm_{m+1}^{(j)}+\ck^{(r,s)}\sum_{j=s}^{m-1}1\otimes \lm_{m+1}^{(j)}
\end{eqnarray*}
Indeed, let $\un[\ldots|x|y|\ldots]$ ($x$ is at the $i$th spot) be one of the summands that enter the expression for $\ck^{(r,s)}(a'_0[\ldots]\otimes a''_0[\ldots])$. Upon applying $\lm^{(i)}$ this summand will not cancel out with anything else if and only if $x$ and $y$ cannot be switched, i.e. in one of the following three cases: 
\[
x=a'_r, \,\, y=a''_s \,\,\,\, \text{or} \,\,\,\, x=a'_k, \,\, y=a'_{k+1} \,\,\,\, \text{or} \,\,\,\, x=a''_{l}, \,\, y=a''_{l+1}
\]
The corresponding summands will contribute to the first, the second, and the third lines in the right-hand side of the formula, respectively.

According to the above formula
\begin{eqnarray*}
\bm^\dagger\cH & = & \sum_{i=1}^{n}\sum_{r=0}^{n-i}\sum_{s=0}^{m}\bm^\dagger\ck^{(r,s)}(\tau_{n+1}^{-r}\ld_{n+1}^{(i)}\otimes \tau_{m+1}^{-s})\\
&& +\sum_{i=1}^{m}\sum_{r=0}^{n}\sum_{s=0}^{m-i}\bm^\dagger\ck^{(r,s)}(\tau_{n+1}^{-r}\otimes \tau_{m+1}^{-s}\ld_{m+1}^{(i)})\\
 & = & \sum_{i=1}^{n}\sum_{r=0}^{n-i}\sum_{s=0}^{m}\lm^{(r+s+1)}\ck^{(r,s)}_0(\tau_{n+1}^{-r}\ld_{n+1}^{(i)}\otimes \tau_{m+1}^{-s})\\
&& +\sum_{i=1}^{n}\sum_{r=0}^{n-i}\sum_{s=0}^{m}\sum_{j=0}^{r-1}\ck^{(r-1,s)}(\lm_{n+1}^{(j)}\tau_{n+1}^{-r}\ld_{n+1}^{(i)}\otimes \tau_{m+1}^{-s})\\
&& +\sum_{i=1}^{n}\sum_{r=0}^{n-i}\sum_{s=0}^{m}\sum_{j=r}^{n-1}\ck^{(r,s)}(\lm_{n+1}^{(j)}\tau_{n+1}^{-r}\ld_{n+1}^{(i)}\otimes \tau_{m+1}^{-s})\\
&& -\sum_{i=1}^{n}\sum_{r=0}^{n-i}\sum_{s=0}^{m}\sum_{j=0}^{s-1}\ck^{(r,s-1)}(\tau_{n+1}^{-r}\ld_{n+1}^{(i)}\otimes \lm_{m+1}^{(j)}\tau_{m+1}^{-s})\\
&& -\sum_{i=1}^{n}\sum_{r=0}^{n-i}\sum_{s=0}^{m}\sum_{j=s}^{m-1}\ck^{(r,s)}(\tau_{n+1}^{-r}\ld_{n+1}^{(i)}\otimes \lm_{m+1}^{(j)}\tau_{m+1}^{-s})\\
&& +\sum_{i=1}^{m}\sum_{r=0}^{n}\sum_{s=0}^{m-i}\lm^{(r+s+1)}\ck^{(r,s)}_0(\tau_{n+1}^{-r}\otimes \tau_{m+1}^{-s}\ld_{m+1}^{(i)})\\
&& +\sum_{i=1}^{m}\sum_{r=0}^{n}\sum_{s=0}^{m-i}\sum_{j=0}^{r-1}\ck^{(r-1,s)}(\lm_{n+1}^{(j)}\tau_{n+1}^{-r}\otimes \tau_{m+1}^{-s}\ld_{m+1}^{(i)})\\
&& +\sum_{i=1}^{m}\sum_{r=0}^{n}\sum_{s=0}^{m-i}\sum_{j=r}^{n-1}\ck^{(r,s)}(\lm_{n+1}^{(j)}\tau_{n+1}^{-r}\otimes \tau_{m+1}^{-s}\ld_{m+1}^{(i)})\\
&& +\sum_{i=1}^{m}\sum_{r=0}^{n}\sum_{s=0}^{m-i}\sum_{j=0}^{s-1}\ck^{(r,s-1)}(\tau_{n+1}^{-r}\otimes \lm_{m+1}^{(j)}\tau_{m+1}^{-s}\ld_{m+1}^{(i)})\\
&& +\sum_{i=1}^{m}\sum_{r=0}^{n}\sum_{s=0}^{m-i}\sum_{j=s}^{m-1}\ck^{(r,s)}(\tau_{n+1}^{-r}\otimes \lm_{m+1}^{(j)}\tau_{m+1}^{-s}\ld_{m+1}^{(i)})
\end{eqnarray*}
Using the formulas
\begin{equation*}
\lm^{(j)}_{n+1}\tau^{-r}_{n+1}=\begin{cases} \tau^{1-r}_{n}\lm^{(n+1+j-r)}_{n+1} & r>j\\ \tau^{-r}_{n}\lm^{(j-r)}_{n+1} & r\leq j \end{cases}
\end{equation*}
\begin{equation}\label{md}
\lm^{(j)}_{n+1}\ld^{(i)}_{n+1}=
\begin{cases}
-\ld^{(i-1)}_{n}\lm^{(j)}_{n+1} \quad 0\leq j< i-1\leq n-1\\
-\ld^{(i)}_{n}\lm^{(j)}_{n+1} \quad 0\leq i< j\leq n-1\\
-\ld^{(i)}_{n}\lm^{(n)}_{n+1} \quad 1\leq i\leq n-1, j=n\\
-\lm^{(i)}_{n+1}\ld^{(i+1)}_{n+1}-\ld^{(i)}_{n}\lm^{(i)}_{n+1} \quad 0\leq i=j\leq n
\end{cases}
\end{equation}
(cf. \cite[Appendix B]{Shk}), and
\begin{equation}\label{Kvsk}
\cK=\sum_{r=0}^{n}\sum_{s=0}^{m}\ck^{(r,s)}(\tau_{n+1}^{-r}\otimes \tau_{m+1}^{-s})
\end{equation}
one shows that 
\begin{eqnarray*}
\bm^\dagger\cH&=&-\cH(\bm\otimes1+1\otimes\bm)-\cK(\U\otimes1+1\otimes\U)\\
&&+\sum_{i=1}^{n}\sum_{r=0}^{n-i}\sum_{s=0}^{m}\lm^{(r+s+1)}\ck^{(r,s)}_0(\tau_{n+1}^{-r}\ld_{n+1}^{(i)}\otimes \tau_{m+1}^{-s})\\
&&+\sum_{i=1}^{m}\sum_{r=0}^{n}\sum_{s=0}^{m-i}\lm^{(r+s+1)}\ck^{(r,s)}_0(\tau_{n+1}^{-r}\otimes \tau_{m+1}^{-s}\ld_{m+1}^{(i)})
\end{eqnarray*}
Thus, to complete the proof of Lemma \ref{l1ap} it remains to prove
that
\begin{eqnarray*}
\V\ck&=&\sum_{r=0}^{n-1}\sum_{s=0}^{m}\sum_{i=1}^{n-r}\lm^{(r+s+1)}\ck^{(r,s)}_0(\tau_{n+1}^{-r}\ld_{n+1}^{(i)}\otimes \tau_{m+1}^{-s})\nonumber\\
&&+\sum_{r=0}^{n}\sum_{s=0}^{m-1}\sum_{i=1}^{m-s}\lm^{(r+s+1)}\ck^{(r,s)}_0(\tau_{n+1}^{-r}\otimes \tau_{m+1}^{-s}\ld_{m+1}^{(i)})
\end{eqnarray*}
(we have changed the order of summation in the original expression).

Set
$
\tilde{\ck}^{(r,s)}_0:=\lm^{(r+s+1)}\ck^{(r,s)}_0(\tau_{n+1}^{-r}\otimes \tau_{m+1}^{-s}).
$
That is, modulo signs,
\begin{multline*}
\tilde{\ck}^{(r,s)}_0(a'_0[a'_1|\ldots |a'_n]\otimes
a''_0[a''_1|\ldots|a''_m])\\
=\pm\un[\underbrace{a'_{n-r+1}|\ldots|a'_{n}|a''_{m-s+1}|\ldots|a''_{m}}_{\text{sum}\,
\text{over}\,(r,s)-\text{shuffles}}|a'_0a''_0|\underbrace{a'_{1}|\ldots|a'_{n-r}|a''_{1}|\ldots|a''_{m-s}}_{\text{sum}\,
\text{over}\,(n-r,m-s)-\text{shuffles}}]
\end{multline*}
It is easy to see that
\begin{multline}\label{aux1}
\sum_{r=0}^{n-1}\sum_{s=0}^{m}\sum_{i=1}^{n-r}\tilde{\ck}^{(r,s)}_0(\ld_{n+1}^{(i)}\otimes 1)
+\sum_{r=0}^{n}\sum_{s=0}^{m-1}\sum_{i=1}^{m-s}\tilde{\ck}^{(r,s)}_0(1\otimes \ld_{m+1}^{(i)})\\
=-\sum_{r=0}^{n}\sum_{s=0}^{m}\sum_{i=r+s+2}^{n+m+1}\ld_{n+m+2}^{(i)}\tilde{\ck}^{(r,s)}_0
\end{multline}
and
\begin{equation*}
-\tau_{n+m+1}^{-j}\ck=\sum_{r+s=j,\,r\leq n,\,s\leq m}\lm^{(0)}_{n+m+2}\tilde{\ck}^{(r,s)}_0 \quad (0\leq j\leq n+m)
\end{equation*}
This, together with (\ref{md}) and the fact that $h_{n+m+1}\lm^{(0)}_{n+m+2}(x)=x$ provided $x$ is of the form $\un[\ldots]$, implies
\begin{eqnarray*}
\V\ck&=&-h_{n+m+1}\sum_{i=1}^{n+m}\sum_{j=0}^{n+m-i}\tau_{n+m+1}^{-j}\ld_{n+m+1}^{(i)}\ck\\
&=&-h_{n+m+1}\sum_{i=1}^{n+m}\sum_{j=0}^{n+m-i}\ld_{n+m+1}^{(i+j)}\tau_{n+m+1}^{-j}\ck\\
&=&h_{n+m+1}\sum_{i=1}^{n+m}\sum_{j=0}^{n+m-i}\sum_{r+s=j,\,r\leq n,\,s\leq m}\ld_{n+m+1}^{(i+j)}\lm^{(0)}_{n+m+2}\tilde{\ck}^{(r,s)}_0\\
&=&-h_{n+m+1}\sum_{i=1}^{n+m}\sum_{j=0}^{n+m-i}\sum_{r+s=j,\,r\leq n,\,s\leq m}\lm^{(0)}_{n+m+2}\ld_{n+m+2}^{(i+j+1)}\tilde{\ck}^{(r,s)}_0\\
&=&-\sum_{i=1}^{n+m}\sum_{j=0}^{n+m-i}\sum_{r+s=j,\,r\leq n,\,s\leq m}\ld_{n+m+2}^{(i+j+1)}\tilde{\ck}^{(r,s)}_0
\end{eqnarray*}
It remains to compare the latter expression to the right-hand side of (\ref{aux1}). It is enough to show that
\begin{eqnarray*}
\sum_{r=0}^{n}\sum_{s=0}^{m}\sum_{i=r+s+2}^{n+m+1}A(i,r,s)=\sum_{i=1}^{n+m}\sum_{j=0}^{n+m-i}\sum_{r+s=j,\,r\leq n,\,s\leq m}A(i+j+1,r,s)
\end{eqnarray*}
for abstract variables $A(x,y,z)$. Changing the order of summation in the right-hand side and shifting the first argument transforms this identity into the following obvious one:
\begin{eqnarray*}
\sum_{r=0}^{n}\sum_{s=0}^{m}\sum_{i=r+s+2}^{n+m+1}A(i,r,s)=\sum_{j=0}^{n+m-1}\sum_{r+s=j,\,r\leq n,\,s\leq m}\sum_{i=j+2}^{n+m+1}A(i,r,s)
\end{eqnarray*}
\hfill $\blacksquare$

\begin{lemma}\label{l2ap} 
\[
\U\cK-\ck(\V\otimes1+1\otimes\V)=
\lm^{(0)}\cH+\lm^{(\ast)}\cH+\ck(\bd\otimes B) 
\] where $\lm^{(\ast)}_{n+1}:=\lm^{(n)}_{n+1}$.
\end{lemma}
\noindent{\bf Proof.} The expression for $\lm^{(0)}\ck^{(r,s)}(a'_0[a'_1|\ldots |a'_n]\otimes
a''_0[a''_1|\ldots|a''_m])$ involves tensors of two types, namely, the ones that begin with $a'_0$ and the ones that begin with $a''_0$. Let us denote the corresponding operators by $\lm^{(0)}\ck^{(r,s)\prime}$ and $\lm^{(0)}\ck^{(r,s)\prime\prime}$; that is, 
\[
\lm^{(0)}\ck^{(r,s)}=\lm^{(0)}\ck^{(r,s)\prime}+\lm^{(0)}\ck^{(r,s)\prime\prime}
\]
Similarly,
\[
\lm^{(\ast)}\ck^{(r,s)}=\lm^{(\ast)}\ck^{(r,s)\prime}+\lm^{(\ast)}\ck^{(r,s)\prime\prime}
\]
with $\lm^{(\ast)}\ck^{(r,s)\prime},\lm^{(\ast)}\ck^{(r,s)\prime\prime}$ having the similar meaning. The key observation is that
\begin{eqnarray*}
\lm^{(\ast)}\ck^{(r,s)\prime}&=&-\lm^{(0)}\ck^{(r+1,s)\prime}(\tau_{n+1}^{-1}\otimes1) \quad (r\leq n-1)\\
\lm^{(\ast)}\ck^{(r,s)\prime\prime}&=&-\lm^{(0)}\ck^{(r,s+1)\prime\prime}(1\otimes\tau_{m+1}^{-1}) \quad (s\leq m-1)\\
\lm^{(\ast)}\ck^{(n,s)\prime}&=&0,\quad \lm^{(0)}\ck^{(r,0)\prime\prime}=0
\end{eqnarray*}
The formulas imply that for $i=1,\ldots, n$ 
\begin{eqnarray*}
\sum_{r=0}^{n-i}\sum_{s=0}^{m}\lm^{(0)}\ck^{(r,s)}(\tau_{n+1}^{-r}\otimes \tau_{m+1}^{-s})
&+&\sum_{r=0}^{n-i}\sum_{s=0}^{m}\lm^{(\ast)}\ck^{(r,s)}(\tau_{n+1}^{-r}\otimes\tau_{m+1}^{-s})\\
=\sum_{s=0}^{m}\lm^{(0)}\ck^{(0,s)\prime}(1\otimes \tau_{m+1}^{-s})
&-&\sum_{s=0}^{m}\lm^{(0)}\ck^{(n-i+1,s)\prime}(\tau_{n+1}^{i}\otimes\tau_{m+1}^{-s})\\
&+&\sum_{r=0}^{n-i}\lm^{(\ast)}\ck^{(r,m)\prime\prime}(\tau_{n+1}^{-r}\otimes\tau_{m+1}^{-m})
\end{eqnarray*}
and for $i=1,\ldots, m$
\begin{eqnarray*}
\sum_{r=0}^{n}\sum_{s=0}^{m-i}\lm^{(0)}\ck^{(r,s)}(\tau_{n+1}^{-r}\otimes \tau_{m+1}^{-s}) &+&\sum_{r=0}^{n}\sum_{s=0}^{m-i}\lm^{(\ast)}\ck^{(r,s)}(\tau_{n+1}^{-r}\otimes \tau_{m+1}^{-s})\\
=\sum_{s=0}^{m-i}\lm^{(0)}\ck^{(0,s)\prime}(1\otimes \tau_{m+1}^{-s})
&-&\sum_{r=0}^{n}\lm^{(0)}\ck^{(r,m-i+1)\prime\prime}(\tau_{n+1}^{-r}\otimes \tau_{m+1}^{i})
\end{eqnarray*}
Using these equalities, it is easy to show that $\lm^{(0)}\cH+\lm^{(\ast)}\cH$ equals
\begin{eqnarray*}
\sum_{i=1}^{n}\sum_{r=0}^{n-i}\lm^{(\ast)}\ck^{(r,m)\prime\prime}(\tau_{n+1}^{-r}\ld_{n+1}^{(i)}\otimes\tau_{m+1})
&+&\sum_{i=1}^{m}\sum_{s=0}^{m-i}\lm^{(0)}\ck^{(0,s)\prime}(1\otimes \tau_{m+1}^{-s}\ld_{m+1}^{(i)})\\
-\sum_{i=1}^{n}\sum_{s=0}^{m}\lm^{(0)}\ck^{(i,s)\prime}(\ld_{n+1}^{(0)}\tau_{n+1}^{-i}\otimes\tau_{m+1}^{-s})
&-&\sum_{i=1}^{m}\sum_{r=0}^{n}\lm^{(0)}\ck^{(r,i)\prime\prime}(\tau_{n+1}^{-r}\otimes \ld_{m+1}^{(0)}\tau_{m+1}^{-i})\\
&+&\sum_{i=1}^{n}\sum_{s=0}^{m}\lm^{(0)}\ck^{(0,s)\prime}(\ld_{n+1}^{(i)}\otimes \tau_{m+1}^{-s})
\end{eqnarray*}
On the other hand, by (\ref{Kvsk}) and the definition of $\ck^{(r,s)\prime}$ and $\ck^{(r,s)\prime\prime}$, $-\U\cK$ equals 
\begin{eqnarray*}
-\sum_{r=0}^{n}\sum_{s=0}^{m}\ld^{(0)}\lm^{(0)}\ck^{(r,s)}(\tau_{n+1}^{-r}\otimes \tau_{m+1}^{-s})&&\\
=-\sum_{r=0}^{n}\sum_{s=0}^{m}\ld^{(0)}\lm^{(0)}\ck^{(r,s)\prime}(\tau_{n+1}^{-r}\otimes \tau_{m+1}^{-s})
&-&\sum_{r=0}^{n}\sum_{s=1}^{m}\ld^{(0)}\lm^{(0)}\ck^{(r,s)\prime\prime}(\tau_{n+1}^{-r}\otimes \tau_{m+1}^{-s})\\
=\sum_{r=0}^{n}\sum_{s=0}^{m}\lm^{(0)}\ck^{(r,s)\prime}(\ld^{(0)}_{n+1}\tau_{n+1}^{-r}\otimes \tau_{m+1}^{-s})
&+&\sum_{r=0}^{n}\sum_{s=1}^{m}\lm^{(0)}\ck^{(r,s)\prime\prime}(\tau_{n+1}^{-r}\otimes \ld^{(0)}_{m+1}\tau_{m+1}^{-s})\\
\end{eqnarray*}
Therefore, we get the following expression for $\lm^{(0)}\cH+\lm^{(\ast)}\cH-\U\cK$
\begin{eqnarray*}
\sum_{i=1}^{n}\sum_{r=0}^{n-i}\lm^{(\ast)}\ck^{(r,m)\prime\prime}(\tau_{n+1}^{-r}\ld_{n+1}^{(i)}\otimes\tau_{m+1})
&+&\sum_{i=1}^{m}\sum_{s=0}^{m-i}\lm^{(0)}\ck^{(0,s)\prime}(1\otimes \tau_{m+1}^{-s}\ld_{m+1}^{(i)})\\
&+&\sum_{s=0}^{m}\lm^{(0)}\ck^{(0,s)\prime}(\bd_{n+1}\otimes \tau_{m+1}^{-s})
\end{eqnarray*}
To finish the proof, it remains to observe that for all $r,s$
\begin{eqnarray*}
\lm^{(\ast)}\ck^{(r,m)\prime\prime}(1\otimes\tau_{m+1})=\ck(h_{n+1}\otimes1),\quad \lm^{(0)}\ck^{(0,s)\prime}=\ck(1\otimes h_{m+1})
\end{eqnarray*}
and, thus, the previous expression equals 
\begin{eqnarray*}
\sum_{i=1}^{n}\sum_{r=0}^{n-i}\ck(h_{n+1}\tau_{n+1}^{-r}\ld_{n+1}^{(i)}\otimes1)
&+&\sum_{i=1}^{m}\sum_{s=0}^{m-i}\ck(1\otimes h_{m+1}\tau_{m+1}^{-s}\ld_{m+1}^{(i)})\\
&-&\sum_{s=0}^{m}\ck(\bd_{n+1}\otimes h_{m+1}\tau_{m+1}^{-s})\\
&=&-\ck(\V\otimes1+1\otimes\V)-\ck(\bd\otimes B) 
\end{eqnarray*}
\hfill $\blacksquare$

Obviously, Lemma \ref{l1} follows from Lemmas \ref{l1ap} and \ref{l2ap}.

\end{document}